\documentstyle[11pt,amssymb]{article}

\newcommand{\id}{{\mbox{Id}}}
\newcommand{\q}{{\bold q}}
\newcommand{\Vect}{{Vect}_{ss,\circ}(\E)}
\newcommand{\M}{{\cal M}_q}
\newcommand{\ad}{{\mbox{Ad}}}
\newcommand{\R}{{\cal R}}
\newcommand{\g}{{\frak g}}

\newcommand{\gl}{{\frak g\frak l }}

\newcommand{\E}{{\cal E}}
\newcommand{\nab}{\noindent}
\newcommand{\C}{{\Bbb C}}
\newcommand{\Z}{{\Bbb Z}}
\newcommand{\Q}{{\Bbb Q}}

\newcommand{\inv}{^{-1}}
\newcommand{\iso}{{\stackrel{\sim}{\longrightarrow}}}

\begin{document}

\setlength{\parskip}{3pt plus 5pt minus 0pt}
\hfill {\small alg-geom/9607008}\break

\centerline{\Large {\bf CONJUGACY CLASSES IN LOOP GROUPS AND}}
\vskip 1mm
\centerline{\Large {\bf $\bold G$-BUNDLES ON ELLIPTIC CURVES}}
\vskip 5mm
\centerline{\sc VLADIMIR BARANOVSKY AND VICTOR GINZBURG}
\vskip 1cm
{\bf 1. Introduction}

Let $\C[[z]]$ be the ring of formal power series and
$\C((z))$ the field of formal Laurent power series, the field of fractions of
$\C[[z]]$. Given a complex algebraic  group $G$, we will 
write $G((z))$ for the group of $\C((z))$-rational points of $G$, thought
of as a formal ``loop group'', and $a(z)$ for an element of $G((z))$.
Let $q$ be a fixed non-zero complex number.
 Define a ``twisted'' conjugation action of $G((z))$ on itself
by the formula
$$ g(z) : a(z)\,\mapsto\, {^g}\!a= g(q\cdot z)\cdot a(z)\cdot g(z)\inv\,. \eqno(1.1)$$
We are concerned with the problem of classifying
the orbits of the twisted conjugation action. If $q=1$  twisted conjugation
becomes the ordinary conjugation, and the problem reduces to
the classification of  conjugacy classes in $G((z))$.

In this paper we will be interested in the case  $|q|<1$. 
 Let $G[[z]] \subset G((z))$
be the subgroup of $\C[[z]]$-points of $G$. A twisted conjugacy
class in $G((z))$ is called {\it integral} if it contains an element
of $G[[z]]$.

Introduce
the elliptic curve $\E=\C^*/q^\Z$.
Our main result is the following
\vskip 2pt

{\sc Theorem 1.2.} {\it Let $G$ be a complex connected 
semisimple algebraic group. Then there is a natural bijection
between the set of integral  twisted conjugacy
classes in $G((z))$ and the set of
isomorphism classes of semi-stable holomorphic principal $G$-bundles on $\E$.}
\vskip 2pt

%
\vskip 2pt

The main reason we are interested in this result is that twisted conjugacy classes
in $G((z))$ may be interpreted as ordinary conjugacy classes in a larger group.
Specifically, the group $\C^*$ acts on $\C((z))$ by field automorphisms rescaling
the variable $z$, i.e., $t\in \C^*$ acts by
$a(z) \mapsto a(t\cdot z)$. This gives a $\C^*$-action on the group
$G((z))$ called ``{\it rotation of the loop}''. Write $\C^*\ltimes G((z))$
for the corresponding semidirect product. Then, 
for any $(q, a)\,,\,(1,g) \in \C^*\ltimes G((z))\,,$ we have
$\,(1,g)\cdot(q, a)\cdot(1,g)\inv=(q,\,^g\!a)\,.$ Thus,
twisted conjugacy classes in $G((z))$ is essentially the same thing as ordinary
conjugacy classes in a Kac-Moody group (=standard central extension
of $\C^*\ltimes G((z))$).

We arrived at theorem 1.2  while trying to find an algebraic version of
the following unpublished analytic result due to Looijenga (cf. [EFK]).
Let $G$ be a connected complex Lie group,  $G(\C^*)_{hol}$ the group of all 
holomorphic maps $a: \C^* \mapsto G$ (possibly with an essential
singularity at $z=0$), and let
$q$ be a fixed non-zero complex number such that $|q|<1$. 
Then Looijenga showed that

{\sc Proposition 1.3.} {\it There is a natural bijection between
the set of {\sl all} twisted-congugacy classes in $G(\C^*)_{hol}$ and the set of isomorphism
classes of {\sl arbitrary} holomorphic $G$-bundles on $\E$.}

\nab
{\sc Proof.} Observe that the pull-back
via the projection $\pi : \C^* \to \C^*/q^\Z=\E$ establishes an equivalence between
the category of $G$-bundles on $\E$ and the category of $q^\Z$-equivariant
holomorphic $G$-bundles on $\C^*$. We associate to 
 $a\in G(\C^*)_{hol}$ the trivial holomorphic $G$-bundle $\C^*\times G \to \C^*\,$
on $\C^*$
with $q^\Z$-equivariant structure given by the action
$\,q :  (z, g) \mapsto (q\cdot z\,,\,a(z)\cdot g)\,.$
The corresponding $G$-bundle on $\E$ will be referred to as the $G$-bundle with 
multiplier $a$. It is easy to see that two $G$-bundles  on $\E$ associated to
two different multipliers are isomorphic if and only if the multipliers are twisted conjugate.
Conversely, it is known  that any holomorphic $G$-bundle  on $\C^*$ is trivial.
The action of the element $q$ on such a trivial bundle has to be of the form
$\,q :  (z, g) \mapsto (q\cdot z, a(z)\cdot g)\,,$ where
$a:  \C^* \to G$ is a holomorphic map (changing trivialization
has the effect of replacing $a$ by a twisted conjugate map).
Hence, every $q^\Z$-equivariant holomorphic $G$-bundle on $\C^*$ can be obtained via the
above construction. $\quad\square$

Although motivation for theorem 1.2 came from loop groups, the result
itself is most adequately understood in the framework
of $q$-difference equations. To explain this assume, for simplicity, that $G=GL_n$.
Given $q\in \C^*$ and $a(z)\in GL_n((z))$, we 
consider a  difference
equation
$$ x(q\cdot z) = a(z)\cdot x(z)\,, \eqno(1.4)$$
where $x(z)\in \C^n((z))$ is the unknown
$\C^n$-valued formal power series. It is clear that if $x(z)$ is a solution
to (1.4) and $g(z)\in GL_n((z))$, then ${\tilde x}(z) :=g(z)x(z)\in \C^n((z))$
 is a solution to a similar equation with $a(z)$ being replaced by
${\tilde a}(z)=g(q z)\cdot a(z)\cdot g(z)\inv\,,$ a twisted conjugate
loop. Therefore classification
of equations (1.4) modulo transformations $x(z) \mapsto {\tilde x}(z)$ reduces
to the classification of the twisted cojugacy classes in $GL_n((z))$.

Equation (1.4) 
should be
regarded as a $q$-analogue of the
first order differential equation
$$z\frac{dx}{dz} = a(z)\cdot x(z)\,, \eqno(1.5)$$
and twisted conjugation (1.1) should be
regarded as a $q$-analogue of the
gauge transformation:
$\,a(z)\mapsto g(z)\cdot a(z)\cdot g(z)\inv + z\frac{dg}{dz}
g(z)\inv\,.$
It is well-known that the classification of  equivalence classes
of equations like (1.5) depends in an essential way on
the type of functions $a$ and $g$ one is considering. If one puts himself
into analytic framework, then $a$ and $g$ are taken to be elements
of $\gl_n(\C^*)_{hol}$ and $GL_n(\C^*)_{hol}$, respectively.
It is well-known and easy to prove
that in this case the differential equation is completely
determined (up to equivalence) by the monodromy of its fundamental
solution. Thus, there is a natural bijection between the set of equivalence 
classes of differential equations of type (1.5) and the set of conjugacy
classes in $G$. This is a differential equation analogue of proposition 1.3.

The situation changes drastically if $\gl_n(\C^*)_{hol}$ and $GL_n(\C^*)_{hol}$
are replaced by formal loops $\gl_n((z))$ and $GL_n((z))$,
 respectively. The classical theory
says that for the equation to be determined by its monodromy it should
have {\it regular singularity} at $z=0$. This is a differential analogue of the
``integrality'' condition in theorem $1.2$. Thus, the $G$-bundle in theorem
1.2 should be thought of as a $q$-analogue of the monodromy of a differential
equation.

Classification of $q$-difference equations (1.4) is 
equivalent to the classification of $D_q$-modules, where by
a $D_q$-module we mean  a finite rank $\C((z))$-module $M$  equipped 
with an invertible $\C$-linear operator $\q: M\to M$ such that
\[\q (f(z)\cdot m) = f(q z)\cdot (\q m)\quad,\quad \forall f
\in \C((z))\,,\,m\in M\,.\]
A $D_q$-module is a module
over a smash product of the group algebra of the group $q^\Z$ with
$\C((z))$.
A $D_q$-module $M$ is said to be integral if there exists
a free $\C[[z]]$-submodule $L\subset M$ of maximal rank
such that $\q(L)\subset L$ and $\q\inv(L)\subset L$.
Integral $D_q$-modules form an abelian category,
 ${\cal M}_q$. It is easy to verify
that tensor product over $\C((z))$ makes $\M$ into a 
tensor category. On the other hand, 
tensor product
of vector bundles makes the category of
degree zero
semistable holomorphic vector bundles on $\E$ into a tensor category,
$\Vect$, with $Hom$'s given by vector bundle morphisms.
The following result is a natural strengthenning of theorem 1.2
in the $G=GL$-case.
\vskip 2pt

{\sc Theorem 1.6.} {\it The tensor category $\M$ is equivalent to $\Vect$.}

\vskip 4pt

\nab
{\sc Acknowledgements.} {\footnotesize
We are grateful to R. Bezrukavnikov and
M. Kontsevich for helpful discussions.}

\vskip 7mm
\pagebreak[3]
{\bf 2. From loop groups to $\bold G$-bundles on} $\bold{\cal E}$.

The ring homomorphism $\C[[z]] \to \C\,,\, f\mapsto f(0)\,$
induces, for any algebraic group $H$,
a natural group homomorphism $H[[z]] \to H$.
Let $H_1[[z]]$ denote the kernel of this homomorphism, 
a ``congruence subgroup''. We use the notation  $H[z]$ and $H[z,z\inv]$ for the
groups of $\C[z]$- and $\C[z,z\inv]$-points of $H$,
 respectively. Thus, $H[z]\subset H[[z]]$ and $H[z,z\inv]\subset H((z))\,.$
Elements of $H[z,z\inv]$ will be referred to as {\it polynomial} loops.

From now on we fix a complex connected
semisimple algebraic group $G$ with Lie algebra $\g$,
 and $q\in \C^*$ such that $|q|<1$.

Our proof of theorem 1.2 consists of several steps. We first assign
to an integral element $a \in G((z))$
a $G$-bundle on $\E$.
The naive idea of using $a$ as a multiplier (cf. proof of 
proposition 1.3) can not be applied here directly, for $a$ is only
a {\it formal} looop, hence, does not give  a holomorphic map in general.
To overcome this difficulty, we prove the following result.

{\sc Proposition 2.1.} {\it For any $a \in G[[z]]$, there exists a Borel subgroup
$B\subset G$ with  unipotent radical $U$, such that $a$ is twisted conjugate
to a polynomial loop of the form $a_0\cdot a_1(z)$ where
$a_0\in B$ and $a_1\in U[z]$.}

To prove the proposition we need some preparations.
Recall that for a semisimple element $s \in G$,
the adjoint action of $s$ on $\g$ has a weight space decomposition
$\g$ $=$ $\bigoplus_{\lambda}\, \g_{\lambda}\,$ where $\g_{\lambda}$
is the eigenspace corresponding to an eigenvalue $\lambda \in \C^*$.

Let $a(z)= a_0\cdot a_1(z)\in G[[z]]$, 
where $a_0 \in G$ is a constant loop and  $a_1(z) \in 
G_1[[z]]$. Write $a_0^{ss}\in G$ for the  the semisimple
part in the Jordan decomposition of $a_0$, and let
$\g$ $=$ $\bigoplus_{\lambda}\, \g_{\lambda}$ be the weight
space decomposition with respect to the adjoint action of $a_0^{ss}$.

{\sc Definition.} The element $a(z)= a_0\cdot a_1(z)$
 is called {\it aligned} if it can be written as a product
$a_0 \exp({x_1z})\exp({x_2z^2})\cdot\ldots\,,$ where $x_i \in \g_{q^i}$.

Note that the product above is finite and gives an element of
$G[z]$. Hence any aligned  element is a polynomial loop.

\vskip 1pt {\sc Lemma 2.2.} {\it For any $a\in G[[z]]$, 
 one can find $g \in G_1[[z]]$ such that
${^g}\!a$ is aligned. }

\nab {\sc Proof.}  Following [BV,pp.31, 68], we will construct a sequence of 
elements $x_i \in \g $ and $y_i \in \g_{q^i}$ as follows.
Note that the exponential map gives a bijection
$z\cdot\g[[z]]\; \iso\; G_1[[z]]$. Therefore we can 
write $a$ in the form $a=a_0\exp({a_1z})\exp({a' z^2})\,,$ where
$a_1 \in \g$ and $a' \in \g[[z]]$. 

Since the operator
$ (q \cdot Ad_{a_0^{-1}} - \id) $ is invertible on $\bigoplus_{\lambda \neq q} 
\g_{\lambda}$, there are uniquely determined elements $x_1 \in \bigoplus_{\lambda \neq q} 
\g_{\lambda}$ and  $y_1\in \g_{q}$ such that 
\[
  (q \cdot Ad_{a_0^{-1}} - \id)(x_1) + a_1 = y_1\,.
\]

We next find $y_2$. To that end, set $g_1=\exp({x_1z})$. 
Then the above equation implies 
 $^{g_1}\!a = a_0\exp({y_1z})\exp(a_2z^2)\exp(a'z^3)\,,$
where $a_2 \in \g$ and $a' \in \g[[z]]$.
Hence there exist uniquely determined
 elements $x_2 \in \bigoplus_{\lambda \neq q^2} 
\g_{\lambda}$ and $y_2 \in \g_{q^2}$ such that
$$
(q^2 \cdot Ad_{a_0^{-1}} - \id)(x_2) + a_2 = y_2\,.
$$
Set $g_2= \exp({x_2}z^2)\exp({x_1z})$. Then the above equation insures
that $^{g_2}\!a=a_0\exp({y_1z})\exp({y_2}z^2)\exp(a_3z^3)\exp(a''z^4)\,,$
where $a_3\in \g$ and $a'' \in \g[[z]]$. Iterating this process we construct
a sequence $\{x_i \in \g\,,\, i=1,2,\ldots\},$
such that setting $g_k:=\exp({x_kz^k})\exp({x_{k-1}z^{k-1}}) \ldots$
$\exp({x_1z})$ we get 
$$
^{g_k}\!a = a_0 \cdot \exp({y_1z})\exp({y_2z^2}) \ldots 
\exp({y_kz^k})\exp({y z^{k+1}})\,,
\eqno(2.3)
$$ 
where $y_i \in \g_{q^i} $ and $y \in \g[[z]]$. Then the 
product $g:=\lim g_k =\ldots\exp({x_kz^k})\cdot$
$\exp({x_{k-1}z^{k-1}})\cdot \ldots\cdot
\exp({x_1z})$ stabilizes since $\g_{q^k}=0$ for all $k>>0$.
Equation (2.3) shows that ${^g}\!a$ is aligned.
$\quad\square$
\vskip 3pt

\nab
{\sc Proof of proposition 2.1.} Choose a maximal torus $T\subset G$ 
 containing $a_0^{ss}$. Let $R \subset Hom(T,\C^*)\,$ be the
set of
roots of $(G,T)$. The subset consisting of the roots $\gamma \in R$
such that 
$\,|\gamma (a_0^{ss})| \leq 1\,$ defines a parabolic $P\subset G$.
Then, for any $i>0$, the subspace $\g_{q^i}$ is contained
in the nilradical of Lie$P$, for $|q|<1$.

Further, we may choose a Borel subgroup $B\subset P$ that contains the unipotent part of 
$a_0$. Let $U$ denote the unipotent radical of $B$.
Then the element
$\exp({y_1z})\exp({y_2z^2}) \ldots \exp({y_kz^k})$
 constructed in the proof of lemma 2.2 belongs to $U[z]$, and the proposition follows.
$\quad\square$
\vskip 3pt

{\sc Lemma 2.5.}{\it
  Let $B$ and  $\widetilde{B}$ be two Borel subgroups with
unipotent radicals $U$, $\widetilde{U}$.
Let  $a=a_0\cdot a_1\,,\,( a_0\in B \,,\,a_1\in U_1[z])\,,$ and 
$\tilde{a}=\tilde{a}_0\cdot \tilde{a}_1\,,\,(
 \tilde{a}_0\in B \,,\,\tilde{a}_1\in U_1[z])\,.$ Then any 
element $g\in G((z))$ such that ${^g}\!a=\tilde{a}$ is  a 
polynomial loop, i.e., $g \in G[z,z\inv]$.}

\nab {\sc Proof.} Multiplying $g$ by an element of $G$ we may assume that 
$B = \widetilde{B}$.
Further, we find a faithful rational representation $\rho: G\to SL_n(\C)$
such that $B$ is the inverse image of the subgroup of upper triangular
matrices in $SL_n(\C)$. Thus, applying $\rho$, we are reduced to proving the
lemma in the case $G=SL_n(\C)$ and $B=\,$upper triangular
matrices.
Thus, from now on, $a$ and $\tilde{a}$ are are assumed to be
 upper  triangular polynomial matrices.
Set $M= max (\deg a, \deg \tilde{a})\,,$ the maximum of the degrees of the corresponding
matrix-valued polynomials. Note that, by assumption, the diagonal entries $a_{ii}$ and 
$\tilde{a}_{ii}$ of the matrices $a$ and $\tilde{a}$ are independent of $z$.

Let
${^g}\!a=\tilde{a}$. We can write $g(z)=$
$\sum_{k\geq k_0} g(k) z^k\,,$ where $g(k)$ are complex $n\times n$-matrices.
Computing the bottom left corner matrix entry  of each side of the
equation $g(qz) a(z) = \tilde{a}(z)g(z)$ yields:
$$
  g(k)_{n,1}\cdot(q^k a_{1,1} - \tilde{a}_{n,n})=0\,.
$$
It follows,  since the diagonal entries of $a, \tilde{a}$ are nonzero,
that there exists at most one $k$, say $k=K_1$, such that 
$g(k)_{n,1} \neq 0$. Using this, we now compute the two matrix entries
standing on $(n-1)\times 1$ and $n\times 2$
places
of each side of the
equation $g(qz) a(z) = \tilde{a}(z)g(z)$. We find that, for any $k >K_1+M$:
$$
  g(k)_{n-1,1}\cdot(q^k a_{1,1} - \tilde{a}_{n-1,n-1})= 0\quad,\quad
g(k)_{n,2}\cdot(q^k a_{2,2} - \tilde{a}_{n,n})=0 \,.
$$
We deduce, as before, that there exists $K_2 >>0$ such that for all $k\geq K_2$, we have
$g(k)_{n-1,1}= g(k)_{n,2}= 0$.

Continuing the process of computing the entries of each side of the
equation $g(qz) a(z) = \tilde{a}(z)g(z)$ along the diagonals
(moving from bottom left corner to top right corner) we prove 
by descending induction on $(i-j)$ that 
$g(k)_{i,j} = 0 $, for all $k\gg0$. $\quad\square$ 

\vskip 4pt
We define a map from integral twisted conjugacy classes in $G((z))$ to
$G$-bundles on the elliptic curve $\E=\C^*/q^\Z$ as follows.
Given an 
element $a\in G((z))$ in an integral twisted conjugacy class,
 we find (proposition 2.1) an aligned element $f\in G((z))$
which is twisted conjugate to $a$. The loop $f$  being polynomial, it gives a well-defined
holomorphic map $f: \C^* \to G$. Hence, we
can associate to $f$ the holomorphic $G$-bundle on $\E$ with multiplier $f$,
see proof of proposition
1.3. If $f'$
is another aligned element which is twisted conjugate to $a$, then by lemma 2.5, $f$ and $f'$
are twisted conjugate to each other via a Laurent polynomial, hence  a holomorphic
 loop. It follows that the $G$-bundles with multipliers $f$ and $f'$
are isomorphic. Thus, we have associated to $a$ a well-defined isomorphism
class of $G$-bundles on $\E$.

{\sc Remark 2.6.} Note that if $a$ is a 
polynomial loop, 
the above map does {\it not} necessarily takes the twisted conjugacy class
of $a$ (in $G((z))$)
to the holomorphic $G$-bundle on $\E$ with multiplier $a$,
even though the latter is defined. We do not claim
that if
$a$ and $a'$ are two 
polynomial loops that are twisted conjugate in $G((z))$ then 
 the $G$-bundles on $\E$ with multipliers $a$ and $a'$ are isomorphic.
Moreover, the following polynomial loops
$$
a={{z\,\,\,\,\,0}\choose {0\,\,z\inv}}
\quad,\quad a'={{z\,\,z\inv}\choose {0\,\,z\inv}}
$$
are twisted conjugate in $SL_2((z))$ by a {\it divergent}
 element, while
the holomorphic $SL_2$-bundles on $\E$ with multipliers
 $a$ and $a'$ are not isomorphic. This is an obstacle for
trying to extend the correspondence of theorem 1.2 beyond
the set of {\it integral} twisted conjugacy classes on the
one hand, and {\it semi-stable} $G$-bundles on the other. 

\vskip 7mm
\pagebreak[3]
{\bf 3. Going to a finite covering.}

Recall that for any positive integer $m$ the field imbedding
$\C((z))\hookrightarrow \C((w))\,,\, z\mapsto w^m\,,$ makes
$\C((w))$ a Galois extension of $\C((z))$ with the Galois group
$\Z/m\Z$. From now on we will write $z^{1/m}$ instead of $w$,
so that $(z^{1/m})^m=z$, and the generator of the Galois group
acts as $\omega: z^{1/m}\mapsto e^{2\pi i/m}z^{1/m}\,.$
Let $G((z^{1/m}))$ denote the group of $\C((z^{1/m}))$-rational points
of $G$. We view $G((z))$ as the subgroup
of $\omega$-fixed points in $G((z^{1/m}))$. We will sometimes write
$a=a(z^{1/m})$ for an element of $G((z^{1/m}))$.

Further, we fix $\tau$ in the upper half-plane, $\mbox{Im} \tau >0$, such that
$q=e^{2\pi i\tau}$. The automorphism $f(z)\mapsto f(q\cdot z)$ of the field
$\C((z))$ can be extended to an automorphism of $\C((z^{1/m}))$ via the
assignment $z^{1/m}\mapsto   e^{2\pi i\cdot\tau/m} z^{1/m}$. This gives rise to
a twisted conjugation action $g: a\mapsto {^g}\!a$ on $G((z^{1/m}))$ 
that extends the one on $G((z))$.

{\sc Definition.}
An element $s \in G$ is said to be {\it reduced} if, for any finite dimensional
rational representation $\rho: G\to GL(V)$, and any 
eigenvalue $\lambda$ of
the operator $\rho(s)$, the equation $\lambda^k=q^l$, (for some $k,l\in \Z\setminus\{0\}$)
implies $\lambda=1,$ i.e., if there are no eigenvalues of the form $\lambda=
e^{2\pi i\cdot\tau\cdot r}\,,\, r\in\Bbb Q\,,$ except $\lambda=1$.

\vskip 2pt
View $G$ as the subgroup of ``constant loops'' in  $G((z^{1/m}))$.

{\sc Theorem 3.1.} {\it Let $a(z) \in G[[z]]$ be an aligned element. Then one 
can find a positive integer $m$ and  $g \in  G((z^{1/m}))$ such that
${^g}\!a$ is a constant loop and moreover the element
${^g}\!a\in G$ is reduced.
}

To prove the theorem, we fix a maximal torus $T\subset G$,
and let $X^*(T)={Hom}_{alg} (T,\C^*)$ and
$X_*(T)={Hom}_{alg} (\C^*,T)\,$ denote the weight and co-weight lattices,
respectively. There is a canonical perfect pairing $\langle\,,\,\rangle:
X^*(T)\times X_*(T)\to \Z\,.$
We first prove

{\sc Lemma 3.2.} {\it For any $s\in T$ there exists ${\phi}\in X_*(T)$
 and an integer
$m\neq 0$ 
such that the following holds:

(i) $s={\phi}(e^{2\pi i\cdot\tau/m})\cdot s_{red}$ where  $s_{red}$ is reduced;

(ii) Let $\alpha\in X^*(T)$. If $\alpha(s)=q^l$ for some $l\in\Z$,
then
$\langle\alpha,{\phi}\rangle/m=l$ and 
$\hphantom{x}\qquad\qquad\alpha(s_{red})= 1.$}

\nab {\sc Proof of Lemma.}  In $\C^*$ consider the subgroup
\[\Gamma=\{ z\in \C^*\;|\; \exists\ k,l\in\Z\;\;\mbox{such that}\;\;z^k=q^l\}\,.\]
 Let $L$ be the subgroup of the weights $\alpha\in X^*(T)$
such that $\alpha(s)\in \Gamma$. Clearly, if $\alpha\in X^*(T)$
and $m\cdot\alpha\in\Gamma$ for some integer $m\neq 0$, then $\alpha\in\Gamma$.
Hence, by the well known structure theorem about subgroups in $\Z^n$ we deduce
that $L$ splits off as a direct summand in $X^*(T)$. Therefore, there is
another lattice $L_{red}\subset X^*(T)$ such that $X^*(T)= L\oplus L_{red}$.
This direct sum decomposition of lattices must be induced by a direct
product decomposition $T=T_1\times T_{red}$, where $T_1$ and
$T_{red}$ are subtori in $T$ such that
$X^*(T_1)=L$ and $X^*(T_{red})=L_{red}$.
 Thus, we have $s=s_1\cdot s'_{red}$, where $s_1\in T_1$
and $s'_{red}\in T_{red}$. 

For any $\alpha\in X^*(T)$, we have by construction $\alpha(s_1)\in \Gamma$.
Furthermore, $\alpha(s)\in \Gamma$ implies $\alpha \in L$, hence $\alpha(T_{red})=1$.
Therefore, for $\alpha\in X^*(T)$ such that $\alpha(s'_{red})\in \Gamma$ we have
$\alpha(s_1\cdot s'_{red})\in \Gamma$, hence $\alpha \in L$, hence
$\alpha(s'_{red})=1$. Thus, $s'_{red}$ is reduced.

View the groups $X^*(T_1)$ and $X_*(T_1)$ as lattices in $
(\mbox{Lie}   T_1)^*$ and
$\mbox{Lie}   T_1$, respectively, so that $X_*(T_1)$ is the kernel
of the exponential map.
Write $s_1= \exp(h) ,$ where $h \in \mbox{Lie}   T_1$. Since
$\alpha(s_1)\in \Gamma$ for any $\alpha\in X^*(T)$, and
elements of $\Gamma$ have the form
$z=e^{2\pi i(\tau\cdot r +r')}\,,$
$r,r'\in \Q\,$, it follows
that $\alpha(h) \in \Q\cdot\tau +\Q\,.$ 
Hence, $h\in \tau\cdot\Q\otimes_{_\Z}
X_*(T_1)+\Q\otimes_{_\Z}
X_*(T_1)\,$. Therefore, there exist ${\phi},{\psi} \in X_*(T_1)$
 and an integer $m$ such that
$h=\frac{\tau}{m}{\phi}+\frac{1}{m}{\psi}\,.$ Thus,
$s_1=exp(h)=\epsilon\cdot{\phi}(e^{2\pi i\cdot \tau/m})\,,$ where
$\epsilon={\psi}(e^{2\pi i/m})$ is an element of order $m$. 
We put $s_{red}=\epsilon\cdot s'_{red}$. Clearly, $s_{red}$ is reduced,
and $s=s_1\cdot s'_{red}= {\phi}(e^{2\pi i\cdot \tau/m})\cdot s_{red}$.

To prove part (ii), 
let $\alpha\in X^*(T)$ be such that $\alpha(s)=q^l$ for some $l\in\Z$.
Then $\alpha\in L$, hence $\alpha(s'_{red})=1$. Furthermore, the equation
$$e^{2\pi i\cdot\tau\cdot l} = q^l =\alpha(s)
=e^{2\pi i\cdot\tau\cdot\langle\alpha,{\phi}\rangle/m+ 2\pi i\cdot\langle\alpha,
{\psi}\rangle/m}$$
yields $\;\tau\cdot(l-\langle\alpha,{\phi}\rangle/m)+\langle\alpha,{\psi}\rangle/m
\in \Z\;.$
It follows, since $\langle\alpha,{\phi}\rangle$ and $\langle\alpha,{\psi}\rangle$ 
are integers,
that $l=\langle\alpha,{\phi}\rangle/m$ and that
$\langle\alpha,{\psi}\rangle/m\in \Z$.
Hence $\alpha(\epsilon)=$
$\alpha({\psi}(e^{2\pi i/m}))=1\,.$
Thus, $\alpha(s_{red})=\alpha(\epsilon)\cdot\alpha(s'_{red})=1\,,$ and
(ii) follows.
$\quad\square$

\nab
{\sc Proof of theorem 3.1.} We choose the Borel subgroup
$B=T\cdot U$ as constructed in the proof of proposition
 2.1.
Thus we have $a(z)= a_0 \exp({x_1 z})\exp({x_2 z^2})\cdot \ldots\cdot
exp({x_k z^k}),$
where $a_0^{ss}\in T$ and $x_i \in \g_{q^i} \subset \mbox{Lie} B\,,$
where $\g_{q^i}$ stands for the $q^i$-eigenspace of $\ad\ a_0^{ss}$.

Applying lemma 3.2 to $s=a_0^{ss}$, we find an integer $m$ and an algebraic
group homomorphism ${\phi}: \C^* \to T$ such that 
$a_0^{ss}={\phi}(e^{2\pi i\cdot\tau/m})
\cdot s_{red}\,.$

For any integer $i\geq 1$ we can write
 $x_i=\sum_\alpha\, x_{\alpha,i}$, where $\alpha$ is a positive
root of
$(G,T)$ such that $\alpha(a_0^{ss})=q^i$ and $x_{\alpha,i}$ is a non-zero
root vector
corresponding to $\alpha$. For such an $\alpha$,
 part (ii) of lemma 3.2 yields
$\alpha({\phi}(z^{1/m}))=z^{\langle\alpha,\phi\rangle/m}= z^i\,.$
We set $g = {\phi}(z^{1/m})\,,$ a well-defined element of the
group $G((z^{1/m}))$. Then, we obtain
\[(\ad\ g)(x_{\alpha,i})= 
\alpha({\phi}(z^{1/m}))\cdot x_{\alpha,i}
 = z^i\cdot x_{\alpha,i}\,.\]
It follows that a similar equation holds for $x_i$ instead of $x_{\alpha,i}$.
From this we deduce 
$$g\inv\cdot\exp(x_iz^i)\cdot g=\exp(x_i)\,.\eqno(3.2.1)$$
Further, let $u$ be the unipotent part of the Jordan decomposition of $a_0$.
Write $u=\exp(y)$. We have $y=\sum_\alpha\, y_\alpha$, where
$y_\alpha\in \mbox{Lie} U$ are root vectors. Since $a_0^{ss}$ commutes with $y$,
we deduce similarly, using lemma 3.2(ii), that $\alpha({\phi})=0$
for any root $\alpha$ such that $y_\alpha\neq 0$. It follows that
$g\inv\cdot u\cdot g=u$.
From this and (3.2.1) we obtain
\[{^g}\!a={\phi}(e^{2\pi i\cdot\tau/m}z^{1/m})\inv\cdot a_0^{ss}\cdot u
\cdot \exp({x_1 z})\exp({x_2 z^2})\cdot \ldots\cdot
exp({x_k z^k})\cdot {\phi}(z^{1/m})\]
\[\quad =
{\phi}(e^{2\pi i\cdot\tau/m})\inv\cdot a_0^{ss}\cdot u
\cdot
\exp({x_1})\cdot\exp({x_2})\ldots
\exp({x_k})\,.\]
Using lemma 3.2(i) we see that ${\phi}(e^{2\pi i\cdot\tau/m})\inv\cdot a_0^{ss}\cdot u=$
$s_{red}\cdot u$. This element is reduced, and the theorem follows.
$\quad\square$

\vskip 1pt 
{\sc Lemma 3.3.}
{\it Let $s\in G$ be reduced. Then any element $g \in G((z))$  such that 
$^g\!s  = s$ is a constant loop.
}

\nab
{\sc Proof.} Consider the  adjoint representation $\rho: G \to GL(\g)$.
We choose a basis in $\g$ such that $\rho(s)$ is an upper-triangular  matrix.
Given $g$ such that $^g\!s = s$, we write $\rho(g)=$
$\sum_{k\geq k_0} g(k) z^k\,,$ where $g(k)$ are complex $n\times n$-matrices.
The same proccess as in the proof of lemma 2.5 gives
equations of the type 
\[g(k)_{m,n}\cdot(q^k s_{n,n} -s_{m,m})=0\quad,\quad k\in \Z\,.\]
Since $s$ is reduced, this implies $g(k)_{m,n}=0$ for all
$k \neq 0$. Hence the image of $g$
in $GL(\g)((z))$ is constant. It follows that $g$ is itself  constant, for
the kernel of the adjoint
representation $G \to GL(\g)$  is finite. $\Box$

{\sc Corollary 3.4.} {\it Let $a\in G[[z]]$ be aligned and $s\in G$ be
reduced. Assume $g\in G((z^{1/m}))$ is such that ${^g}\!a=s$.
Then $g\in G[z^{1/m},z^{-1/m}]$ is a Laurent polynomial loop in
$z^{1/m}$. Furthermore, $\theta=g(e^{2\pi i\cdot\tau/m} z^{1/m})g(z^{1/m})^{-1}$
is a constant loop, and $\theta^m=1$.}

\nab
{\sc Proof.} The first claim follows from lemma 2.5.
To prove the second claim, 
recall the Galois automorphism $\omega: f(z^{1/m})\mapsto
f(e^{2\pi i\cdot\tau/m} z^{1/m})$ on $\C((z^{1/m}))$. We apply the induced
automorphism of $G((z^{1/m}))$ to the equation ${^g}\!a=s$.
The RHS being independent of $z$, and $a$ being fixed by $\omega$,
we get $^{\omega g}\!a=s$. This equation together with the original one,
${^g}\!a=s$, yields $^\theta\!s=s$, where
$\theta=g(e^{2\pi i\cdot\tau/m} z^{1/m})g(z^{1/m})^{-1}\,.$ Hence, $\theta$
 is a constant loop, by lemma 3.3.
Further, applying the automorphism $\omega$ several times
to the first equation below we get
a sequence of equations
\[\theta=(\omega g)\cdot g\inv\enspace,
\enspace\theta=(\omega^2 g)\cdot(\omega g)\inv\;,\ldots,\;
\theta=(\omega^{m} g)\cdot(\omega^{m-1} g)\inv\,.\]
Since $\omega^{m}=\id$, taking the product of  
all these equations yields $\theta^m=1$.
$\quad\square$

\vskip 4pt

We fix two generators, an ``$a$-cycle'' and a ``$b$-cycle'',
 of the fundamental 
group $\pi_1(\cal E)$ as follows. The $a$-cycle is defined to be
the image of a generator of    
$\pi_1 (\C^*)=\Z$
under the imbedding $\pi_1 (\C^*) \hookrightarrow \pi_1(\cal E)$ 
induced by the projection $\C^* \to \C^*/q^\Z=\E\,.$ The $b$-cycle is 
the image of the segment
$[1,q] \subset \C^*$ under the projection.

Given an integer $m\neq 0$, set $^m\!\E= \C^*/q^{\frac{\Z}{m}}$,
and let $^m\!\pi: {^m\!\E}
 \to \E\,,\, z\mapsto z^m\,$ be the natural projection.
Thus, $^m\!\E$ is an elliptic curve and the map $^m\!\pi$ is an $m$-sheeted
Galois covering with the Galois group $\Z/m\Z$ acting as monodromy
around the $a$-cycle.

\vskip 2pt
{\sc Proposition 3.5.} {\it Let $a\in G[z]$ be an aligned element and
$P$ the
principal $G$-bundle on $\E$ with multiplier $a$. Then
there exists a positive integer $m$ such that

(i) The bundle $^m\!\pi^*P$ is isomorphic to the holomorphic
$G$-bundle on $^m\!\E$ with a reduced constant multiplier $s\in G$;

(ii) Let $\tilde \nabla$ be the holomorphic connection on
$^m\!\pi^*P$ transported via the isomorphism (i)
from the trivial connection $d$ on the $G$-bundle with multiplier $s$.
Then $\tilde\nabla$ descends to a well-defined holomorphic 
connection $\nabla$ on $P$.  The latter has finite monodromy
around $a$-cycle and a reduced monodromy around $b$-cycle.}

\nab {\sc Proof.} By theorem 3.1 there exists an element $g
\in G(z^{1/m})$ such that ${^g}\!a=s$ is a constant loop, and
 $s \in G$   is reduced. By corollary 3.4, 
$g=g(z^{1/m})$ is a Laurent polynomial in
$z^{1/m}$. Hence,  $g$ may be viewed as a well-defined $G$-valued
regular function on the $m$-fold covering of $\C^*$. 
Let $P$ be the $G$-bundle on $\E$ with
multiplier $a$.
It follows
that the pull-back, $^m\!\pi^*P$, has a multiplier which is a
twisted conjugate
of $s$. This proves part (i).

To prove (ii), recall that any $G$-bundle with a constant multiplier $s$
has a natural flat holomorphic connection
 which is given (in the trivialization
on $\C^*$ corresponding to $s$) by the deRham differential $d$.
We transport this connection to $^m\!\pi^*P$ via the
isomorphism given by the loop $g$. The connection $\tilde\nabla=g\inv\circ d \circ g$
thus obtained descends to to a connection on $P$ if and only if it is invariant
under the Galois action of $\Z/m\Z$. But by corollary 3.4 we have
$\omega g= \theta\cdot g$, hence we get
\[\omega(\tilde\nabla)=
(\theta\cdot g)\inv\circ d\circ (\theta\cdot g) 
= g\inv\cdot( \theta\inv\circ d\circ\theta)\cdot g =
g\inv\circ d\circ g=\tilde\nabla\,,\]
since $\theta$ commutes with $d$.
Thus, $\tilde\nabla$ is fixed by the Galois action.

To compute the monodromy, note that $g\inv$ is a flat section of the connection
$\tilde\nabla$. Hence the monodromy of $\tilde\nabla$ around  $b$-cycle equals
$g(e^{2\pi i\cdot\tau/m} z^{1/m})\cdot$
$g(z^{1/m})^{-1} =\theta\,.$
Since the covering $^m\!\pi: {^m\!\E}
 \to \E$ has no monodromy around $b$-cycle
and has finite monodromy around $a$-cycle, it follows
that $\nabla$ also has  monodromy  $\theta$ around $b$-cycle
and has finite monodromy around $a$-cycle.$\quad\square$

Given a finite dimensional rational $G$-module $V$,
write $V_{_P}$ for the associated vector bundle on $\E$ corresponding
to a principal $G$-bundle $P$.

\vskip 2pt
{\sc Lemma 3.6.} {\it Let $P$ be the $G$-bundle with an aligned multiplier,
and $\nabla$ the connection on $P$ constructed in proposition 3.5.
Then, for any rational representation $\phi: G\to GL(V)$,
every holomorphic section of the associated vector bundle $V_{_P}$ is flat with respect
to the induced connection on $V_{_P}$.
}

\nab{\sc Proof.} Since the connection $\nabla$
constructed in proposition 3.5
was obtained from a connection $\tilde\nabla$ on $^m\!\pi^*P$,
the claim for $V_{_{P}}$
is equivalent to a similar claim for the vector
bundle $^m\!\pi^*V_{_{P}}$. This vector bundle is isomorphic to the vector bundle $\cal V$
on $^m\!\E$  with multiplier $\phi(s)$, so that the connection $\tilde\nabla$
is isomorphic to the trivial connection $d$. Thus, proving the claim
amounts to showing that any holomorphic section of the vector bundle $\cal V$
with  multiplier $\phi(s)$ is constant.

To that end, write the matrix $\phi(s)$ in 
Jordan form  $\phi(s)= \bigoplus_i J(\lambda_i, n_i)\,,$ where
$J(\lambda_i, n_i)$ is the $(n_i\times n_i)$ Jordan block with 
eigenvalue $\lambda_i$. This gives the corresponding vector
bundle decomposition ${\cal V}={\bigoplus}_{i} {\cal V}_{i}$ where ${\cal V}_i$
is the vector bundle with multiplier $J(\lambda_i, n_i)$.
If $L_i$ denotes the line bundle with multiplier $\lambda_i$,
then there is a canonical vector bundle imbedding $L_i\hookrightarrow
{\cal V}_i$. Furthermore, one can prove (using, e.g., the Fourier-Mukai transform)
that the imbedding induces an isomorphism
$\Gamma(^m\!\E,\,L_i)\,\stackrel{\sim}{\to}\,\Gamma(^m\!\E,\,{\cal V}_i)\,$
of the spaces of global sections. Hence, any holomorphic section
of ${\cal V}_i$ comes from a holomorphic section of $L_i$.
But $L_i$ is a degree zero line bundle, hence has a non-zero section
only if it is the trivial bundle, i.e if $\lambda_i=q^m$.
Observe now that $\lambda_i$ is an eigen-value of the matrix
$\phi(s)$. Since $s\in G$ is reduced, equation $\lambda_i=q^m$
implies $\lambda_i=1$. But then the only holomorphic section
of $L_i$ is a constant section.
The latter is annihilated by the deRham differential $d$, and the lemma is proved.
$\quad\square$

\vskip 2pt
{\sc Proposition 3.7.} {\it Let $a, a_1\in G((z))$ be two aligned
elements. If the $G$-bundle on $\E$ with multiplier
$a$ is isomorphic to the $G$-bundle  on $\E$ with multiplier $a_1$,
then $a$ is twisted conjugate to $a_1$ via a polynomial loop.}

\nab
{\sc Proof.} By theorem 3.1,  there exist an integer $m\geq 1$ and elements
$g, g_1\in G((z^{1/m}))$ such that
$$^ga=s\quad,\quad ^{g_1}a_1=s_1\quad\mbox{where}\enspace s,s_1\in G\enspace
\mbox{are reduced.}\eqno(3.7.1)$$
Let $\nabla$, $\nabla_1$ be the holomorphic connections
on the $G$-bundles on $\E$ with multipliers $a$ and $a_1$, respectively, constructed
in proposition 3.5. The monodromies of the connections around  $a$-cycle,
are equal to $s$ and $s_1$, respectively, and the monodromies 
around  $b$-cycle are are equal to
$\theta$ and $\theta_1$, respectively. By proposition 3.5 we have
$\theta^m=\theta^m_1=1\,.$
If the $G$-bundles with multipliers $a$ and $a_1$
are isomorphic, then we may view $\nabla_1$ as another holomorphic
connection
on the  $G$-bundle $P$ with multiplier $a$.

  Since the cotagent bundle on $\cal E$ is trivial
the difference $X=\nabla_1 - \nabla$ may 
be viewed as a holomorphic section, $X$, of the adjoint
bundle $\g_{_P}$.
Since $s$ is reduced the section $X$
is flat with respect  to $\nabla$, by lemma 3.6. Let
$p: \tilde{\E} \to \E$ be a universal cover of $\E$. The bundle
$p^*P$ on $\tilde{\E}$ has a horisontal holomorphic
section. This section gives a trivialization of $p^*P$ such that,
in the induced trivialization of $p^*\g_{_P}$,
the pull-back $p^*X$ is a constant element $x\in \g$.
Observe that in general, any element $y\in \g\,$ gives rise in this way
to a flat multivalued section of $\g_{_P}$, and 
the monodromy of this section
around $a$- and $b$-cycles is equal to
$\ad \theta (y)$ and $\ad s (y)$, respectively.
It follows, since $X$ is a single-valued flat section of $\g_{_P}$
without monodromy, that $x$ commutes with both $\theta$ and $s$.
Hence, equation $\nabla_1 = \nabla + X$ shows that the monodromy
of the connection $\nabla_1$ is given by the formulas
$$\theta_1= \exp(x)\cdot\theta\quad,\quad
s_1= \exp(\tau x)\cdot s\,.\eqno(3.7.2)$$
From  these formulas and the equations $\theta_1^m=\theta^m=1$ we deduce
$\exp(m\cdot x)=1$. Thus, we may find a maximal torus $T$ containing $\theta, \theta_1$,
and $\phi\in X_*(T)$, such that $x=\phi/m$ (cf. proof of lemma 3.2).

Clearly,
$\phi(z^{1/m})$ is a well defined element of $G((z^{1/m}))$, and from the
first formula in (7.3.2) we deduce $\phi(e^{2\pi i\cdot\tau/m}z^{1/m})\cdot s\cdot
\phi(z^{-1/m}) = s_1\,.$ Recall the notation of (7.3.1), and put $f(z^{1/m})=
g_1(z^{1/m})\inv\cdot \phi(z^{1/m})\cdot g(z^{1/m})\in G((z^{1/m}))\,.$
We claim that $f\in G((z))$. To prove this, it suffices to show that
$f(e^{2\pi i\cdot\tau/m}z^{1/m})= f(z^{1/m})$.
The latter  follows from the chain of equalities:
$$
\begin{array}{l}
f(e^{2\pi i\cdot\tau/m}z^{1/m})=
g_1^{-1}(e^{2\pi i\cdot\tau/m} z^{1/m})\cdot \phi(e^{2\pi i\cdot\tau/m} z^{1/m})
\cdot g (e^{2\pi i\cdot\tau/m} z^{1/m}) =\\
g_1^{-1}(z^{1/m})\cdot  \theta_1^{-1}\cdot  \exp(x)\cdot \phi(z^{1/m})\cdot
\theta\cdot g(z^{1/m})=\\
g_1^{-1}(z^{1/m})\cdot  \theta_1^{-1} \cdot \exp(x)\cdot \theta\cdot \phi(z^{1/m})
\cdot g(z^{1/m}) =\\
 g_1^{-1}(z^{1/m})\cdot \phi(z^{1/m})\cdot g(z^{1/m})=f(z^{1/m})\,.
\end{array}
$$
Finally, using (7.3.1) we calculate
\[^fa\,=\,^{g_1\inv\cdot\phi\cdot g}\!a
\,=\,^{g_1\inv\cdot\phi}\!s\,=\, ^{g_1\inv}\!\!s_1\,=\, a_1\,.\]
Thus, $a$ and $a_1$ are twisted conjugate by an element of $G((z))$.
Lemma 2.5 completes the proof. $\quad\square$

\vskip 7mm
\pagebreak[3]
{\bf 4. Semistable G-bundles and holomorphic connections.}

 Recall that $G$ is a complex connected semisimple group.
 For the definition and properties of semistable holomorphic
$G$-bundles on an elliptic curve we refer to [R] and [RR].

\vskip 2pt
{\sc Proposition 4.1.}
{\it A holomorphic principal $G$-bundle over an elliptic curve is semistable if
and only if it has a holomorphic connection (necessarily flat).}

\nab{\sc Proof.} The ``if'' part is a corollary of the main result of [B].
The ``only if'' part follows from theorem 4.2 below. $\quad\square$

\vskip 2pt
{\sc Theorem 4.2.} {\it For any
  semistable $G$-bundle $P$ on $\E$, there exists
a  holomorphic 
connection on $P$ with finite order monodromy around $a$-cycle
such that, for any rational $G$-module $V$,
every holomorphic section of the associated vector bundle $V_{_P}$ is flat with respect
to the induced connection on~$V_{_P}.$
}
   
\nab
{\sc Proof.} We choose and fix a faithful rational representation
$G\to GL(V)$. By a theorem of Ramanan and Ramanathan
[RR], semistability of $P$ implies semistability of $V_{_P}$. By the classification
of semistable vector bundles on $\E$, due to
Atiyah [A], any semistable vector bundle is isomorphic to the vector bundle
with a constant multiplier. Hence, the bundle $V_{_P}$ has constant multiplier
$a\in GL(V)$. View $a$ as an element of the semisimple group $PGL=PGL(V)$, and let
$P_a$ be the principal $PGL$-bundle with (constant) multiplier $a$.
By construction, the $PGL(V)$-bundle $P_a$ is induced from the $G$-bundle
$P$ via the composition $G\to GL(V)\to PGL(V)\,.$

We may regard the element $a\in PGL$ as a constant aligned loop
in $PGL((z))$. Applying
proposition 3.5, we see that there is an integer $m\geq 1$ and
a reduced element $s\in PGL$ such that the bundle $^m\!\pi^*P$ on
$^m\!\E$ is isomorphic to the
$PGL$-bundle on $^m\!\E$ with multiplier $s\in PGL$.
 Let $\nabla$ be the connection
on $P_a$ constructed in proposition 3.5. 

We claim  that the
connection $\nabla$ on $P_a$ arises from a holomorphic
connection on the $G$-bundle $P$ via the composite homomorphism 
${\rho}: G\rightarrow GL(V)\to PGL(V)\,.$
Note that this composition has finite kernel, so that the induced
canonical map
$i: P\to P_a$ is an immersion. Let $TP_a$ be the tangent bundle
on $P_a$. Our claim is equivalent to saying that
the distribution in $TP_a$ formed by  the horisontal subspaces 
of the connection $\nabla$ is tangent to the immersed submanifold
$i(P)\subset P_a$. Observe that the canonical map $i: P\to P_a$
gives rise to a holomorphic section  $\nu: \E=P/G \to P_a/{\rho}(G)$.
The horisontal distribution is tangent to $i(P)$ if and only
if $\nu$ is a horisontal section.

To show the latter, we apply
Chevalley's theorem [S, Theorem 5.1.3] to the
algebraic subgroup ${\rho}(G)\subset PGL$. The theorem says that
we can find a rational representation $\phi: PGL\to GL(E)$ and a 1-dimensional
subspace ${\bold l} \subset E$ such that ${\rho}(G) = \{ g \in PGL\; | \;
\ \phi (g)(\bold l)= \bold l\}$. Notice, that since $G$ is semisimple, it 
 stabilises a vector $l \in \bold l$. Hence, the assignment $g\mapsto g(l)$
gives rise to an imbedding $PGL/{\rho}(G)\hookrightarrow E$.
Now let $E_{_{P_a}}$ be the associated
vector
bundle corresponding to $E$, equiped with the connectioned induced by $\nabla$.
 The imbedding $PGL/{\rho}(G)\hookrightarrow E$
gives rise to an imbedding $P_a/{\rho}(G) \hookrightarrow E_{_{P_a}}$
compatible with the connections. To show that $\nu$ is horisontal,
it suffices to show that its image under the above imbedding is a flat section.
But this  image is a holomorphic section of  $E_{_{P_a}}$.
By lemma 3.6, any holomorphic section of the vector
bundle $E_{_{P_a}}$ is flat with respect to the connection
on $E_{_{P_a}}$ induced by $\nabla$.  This proves that $\nu$ is horisontal, so that
the horisontal distribution on $TP_a$ is tangent to $i(P)$ and the connection
$\nabla$ comes from a holomorphic $G$-connection on $P$.

Observe further
that the connection $\nabla$ on $P_a$ has finite monodromy around
$a$-cycle. The map $i: P\to P_a$ being an immersion with finite fibers,
it follows that the $G$-connection on $P$  also has finite monodromy around
$a$-cycle.

Finally, it remains to show that there exists a holomorphic connection on $P$
with finite order monodromy around $a$-cycle
such that, for any rational $G$-module $V$,
every holomorphic section of the associated vector bundle $V_{_P}$ is flat with respect
to the induced connection on $V_{_P}$. We do not claim that the connection
we have constructed has this property. Instead we proceed as follows.
We first use the connection that we have constructed above
to prove that any
semistable $G$-bundle on $\E$ is isomorphic to a $G$-bundle with an aligned 
 multiplier. This will be done in the proof of theorem 4.3 below. We can then 
apply proposition 3.5(ii) and lemma 3.6 to get a connection on $P$ with
all the required
properties.$\quad\square$

\vskip 2pt

{\sc Theorem 4.3.} {\it
 A $G$-bundle on $\cal E$ is semistable if and only if it is
isomorphic to the $G$-bundle with an
aligned multiplier $a \in G[z]$.}

\nab
{\sc Proof.} By  proposition 3.5(ii), any $G$-bundle $P$
with an aligned multiplier has a holomorphic connection. Then, the ``if'' part of
proposition 4.1 (due to Biswas) implies that $P$ is semistable.

Conversely, let $P$ be a semistable $G$-bundle. By theorem 4.2, we can
equip $P$ with a holomorphic connection that has monodromies $\theta, b\in G$
around the $a$- and  $b$-cycles respectively, such that $\theta^m=1$ for
some integer $m\geq 1$. Observe that the elements $\theta$ and $b$ commute,
for $\pi_1(\E)$ is an abelian group. Hence, there is a maximal torus
$T\subset G$ such that $\theta, b^{ss}\in T$.
As in the proof of proposition 2.1, we choose a Borel subgroup
$B\supset T$ such that $b\in B$ and $|\alpha(b^{ss})|\leq 1$, for any
positive (with respect to $B$)  root $\alpha$. 

Further, since $\theta^m=1$ there exists ${\phi}\in X_*(T)$ such that
$\theta={\phi}(e^{2\pi i/m})\,.$ Set $g={\phi}(z^{1/m})\inv$,
a well-defined polynomial loop in $G((z^{1/m}))$. We put
$a={^g\!b}\in G((z^{1/m}))$.
 We have $g(e^{2\pi i/m}z^{1/m})=\theta\inv g(z^{1/m})\,.$
Since $\theta$ commutes with $b$, we deduce that
$a(e^{2\pi i/m}z^{1/m})=a(z^{1/m})$. It follows that 
$a$ is fixed by the Galois group, hence, $a\in G((z))$.

Let $U$ be the unipotent radical of $B$. We have $b=b^{ss}\cdot u$ where
$u\in U$. Hence, the condition $|\alpha(b^{ss})|\leq 1$ for any
positive  root $\alpha$, insures that $a=$
$^g\!b=b^{ss}\cdot a_1$
where $a_1\in U_1[[z]]$.
Moreover, since $g$ is a polynomial loop we have $a_1\in U_1[z]$.
By proposition 1.3, the element $a$ is twisted conjugate in $G((z))$
to an aligned element $a'$. Using lemma 2.5 and the fact that $a\in B\cdot U[z]$,
we see that $a$ is twisted conjugate to $a'$ via a polynomial loop. 
Thus, there is an element $f \in G[z^{1/m}, z^{-1/m}]$ such that
\[^fa'= b\quad,\quad f(e^{2\pi i/m}z^{1/m})f(z^{1/m})\inv= \theta\,.\]
These equations show (see proof of proposition 3.5)
that the $G$-bundle $P'$ with multiplier $a'$ has
 a holomorphic connection with the monodromies $\theta$ and $ b\in G$
around $a$- and  $b$-cycle, respectively. 
Thus $P$ and $P'$ are two $G$-bundles with connections that have
the same monodromy. 
Since a holomorphic $G$-bundle with connection is determined, up to isomorphism,
by the monodromy representation, we deduce that
 $P\simeq P'$, and the theorem follows. $\quad\square$

\vskip 2pt
\nab
{\sc Proof of theorem 1.2.} Proposition 3.5 shows that the $G$-bundle
associated to any integral twisted conjugacy class in $G((z))$ via the 
procedure described at the end of $\S 2$ has a holomorphic connection,
hence is semistable, due to proposition 4.1. Theorem 4.3 insures
that the map $\{${\it integral twisted conjugacy classes}$\}
\,\longrightarrow\;\{$ {\it isomorphism classes of 
semistable $G$-bundles} $\}\;$ is surjective.
Injectivity of the map follows from proposition 3.7. $\quad\square$

\vskip 2pt
\nab
{\sc Sketch of proof of theorem 1.6.} We note first that although
$GL_n$ is not a semisimple group, lemma 3.2 and all other
results of sections 3 - 4 hold for $G=GL_n$. In particular, any holomorphic
degree zero semistable vector bundle on $\E$ 
has a holomorphic connection with finite 
order monodromy around $a$-cycle and
a reduced monodromy around $b$-cycle
(since any semistable vector bundle is isomorphic to the
one with a constant multiplier [A],
it has a holomorphic connection  without monodromy
around $a$-cycle at all. M
oreover, it is easy to achieve that the monodromy
matrix around $b$-cycle has no eigenvalues of the
form $q^k\,,\, k\in \Z\setminus\{0\}\,.$ But in order 
to ensure that the monodromy around $b$-cycle
is reduced one still has to allow a finite 
monodromy around  $a$-cycle).
In this way, one shows that an
analogue of theorem 1.2 holds for $G=GL_n$.

Now let $V$ be a rank $n$ and
degree zero
semistable vector bundle on $\E$,
 and $\nabla$ a holomorphic connection
on $V$ with order $m$ monodromy around $a$-cycle.
 Let $\pi: \C^* \to \E$ and $^m\!p: \C^*\to\C^*\,,\,
z\mapsto z^m$ denote the projections.
We know that  
the pull-back $^m\!p^*(\pi^*V)$ can be
trivialized by means of a flat (with respect 
to the pull-back connection) frame.
If $\nabla'$ is 
another holomorphic connection 
on $V$ with finite order monodromy around $a$-cycle, then we can
choose a large enough integer $m$ and two trivializations
of $^m\!p^*(\pi^*V)$ that are flat with respect to
$^m\!p^*\nabla$ and $^m\!p^*\nabla'$, respectively.
We claim that the transition matrix between any two such
trivializations is given by a  holomorphic map $a: 
\C^*\to GL_n$ with a pole 
at $z=0$. To prove the claim, we may choose and fix
 one holomorphic connection 
$\nabla^\circ$ on $V$ that has finite order monodromy 
around $a$-cycle and, moreover,
a reduced monodromy around $b$-cycle. Using an argument
similar to that in the proof of proposition 3.7 one shows
that
the transition matrix between the trivialization corresponding to a 
flat frame with respect to $\nabla^\circ$ and the trivialization corresponding to a 
flat frame  with respect to 
 any other connection with finite monodromy
around $a$-cycle is given by a holomorphic
loop with a pole at $z=0$. This proves the claim.

Next, we define a class of ``moderate'' holomorphic sections of
$\pi^*V$ as follows. Write $\R$ for the ring of holomorphic functions
on $\C^*$ with a pole at $z=0$.
Let $s$ be a holomorphic
section of $\pi^*V$, let $\nabla$ be a holomorphic connection
on $V$ with finite order monodromy around $a$-cycle, and
let $\,s_1,\ldots, s_n\,$ be a frame of flat sections
of $^m\!p^*(\pi^*V)$, for some $m$. 
We say that $s$ is moderate with respect to $\nabla$
if $\,^m\!p^*s\,$ is expressed as a linear combination of the
sections
$\,s_1,\ldots, s_n\,$ with coefficients in $\R$. It is clear that
this definition does not depend on the choice of a flat frame
and the choice of finite covering $^m\!p: \C^*\to\C^*$. Moreover,
the claim proved in the previous paragraph insures that
$s$ is moderate with respect to $\nabla$ if and only if
it is moderate with respect to $\nabla'$ provided
both connections have finite order monodromy around $a$-cycle.
Thus, there is a well defined notion of a {\it moderate}
holomorphic section of $\pi^*V$. Clearly, moderate
sections form an $\R$-module, to be denoted $V(\R)$.

In the previous setting, assume that $V$ is the vector bundle
with a constant multiplier $a\in GL_n(\C)$, 
and $\nabla=d$ is the trivial connection
on $V$.
Then the tautological trivialization of $\pi^*V$ is flat.
Furthermore, 
it is clear that a section of $\pi^*V$ is moderate if and only if 
all its coordinates (relative to the
trivialization) belong to $\R$. It follows that $V(\R)\simeq \R^n$
is
a free $\R$-module. The $q^\Z$-equivariant
structure on $\pi^*V$, see proof of proposition 1.3, provides an operator 
$\q: V(\R)\to V(\R)$. Viewed as an operator $\R^n\to \R^n$ via the 
 isomorphism
$V(\R)\simeq \R^n$, this operator has the form
$$\q: f(z)\mapsto a\inv\cdot f(q\cdot z)\quad\mbox{where}
\quad a=\mbox{ multiplier}\,.\eqno(4.4)
$$

Let $\Vect$ be the tensor category of degree zero holomorphic
semistable vector bundles on $\E$.
We define a functor
$F: \Vect \to \M$ by the assignment
\[V  \mapsto F(V)=\C((z))\otimes_{_\R} V(\R)\,,\]
where the operator $\q:F(V)\to F(V) $ extends the one introduced above.
Formula (4.4) shows that $\q$ and $\q\inv$ preserve the standard lattice
$L=\C[[z]]^n\subset \C((z))^n\,$. Thus, $F(V)\in \M$.

Let $V,V'\in \Vect$.
It is clear from construction that there is a natural map
$V(\R)\otimes_{_\R} V'(\R)\to (V\otimes V')(\R)$. This map is
in effect an isomorphism. To see this, choose trivializations
of $\pi^*V$ and $\pi^*V'$ corresponding to constant multipliers. Then
the tensor product trivialization on $\pi^*(V\otimes V')$ corresponds
to the tensor product of the multipliers. We have seen that in these
trivializations one has $V(\R)\simeq \R^{rk V}\,,$
$V'(\R)\simeq \R^{rk V'}$ and $(V\otimes V')(\R)\simeq \R^{rk V\cdot rk V'}\,.$
It follows that 
the natural map above becomes the standard
isomorphism $\R^{rk V}\otimes_{_\R} \R^{rk V'}\,\iso\,$
$
\R^{rk V\cdot rk V'}\,.$ Applying $\C((z)\times_{_\R}(\bullet)$
we deduce that $F$ is a tensor functor.

Recall now that both $\M$ and $\Vect$ are rigid tensor
categories, cf. [DM], (in particular, $\Vect$ is an abelian
$\C$-category with finite dimensional $Hom$'s). 
To prove that a tensor functor between
two rigid tensor
categories is an equivalence, it suffices to show it is fully faithful.
The analogue of theorem 1.2 for $G=GL_n$ insures that the
 functor $F$ is full. It remains to show that, for any $V',V\in \Vect$,
the functor $F$ induces an isomorphism
 $\,Hom_{_{\Vect}}(V',V)= Hom_{_{\M}}(F(V'), F(V))\,.$ Using the duality
functor on rigid tensor
categories one  reduces to the case $V'= {\bold 1}_{_\E}$, the unit of the
tensor category $\Vect$, i.e. the trivial one-dimensional vector
bundle.
In that case we have
\[\,Hom_{_{\Vect}}({\bold 1}_{_\E}, V)= \Gamma(\E,\,V)= \Gamma_{hol}(\C^*,\,\pi^*V)^\q\,\]
 is the fixpoint  subspace of the $q^\Z$-action on the space 
of all holomorphic sections of $\pi^*V$. On the other hand, we have
$\,Hom(F({\bold 1}_{_\E}), F(V))= F(V)^\q\,,$ is the
 $\q$-fixpoint space  of the operator $\q$. In the trivialization
corresponding to a constant multiplier $a$ the operator $\q$ is given
by formula (4.4). Therefore, we find
\[ F(V)^\q= \{f(z)\in \C((z))^n\;|\; 
f(q\cdot z) = a f(z)\}\]
Expanding $f$ as a formal Laurent series $f(z)=\sum_{k\geq k_0} f_k z^k\,,$
shows that $f\in F(V)^\q$ if and only if, 
 for any $k$, the coefficient
$f_k\in \C^n$ is an eigenvalue of the operator $a$ with eigenvalue $q^k$.
It follows that $f_k$ must vanish for all $k>>0$. Hence $f$ is a Laurent
polynomial
section, $f\in V(\R)^\q$. Similarly, any element of $\Gamma_{hol}(\C^*,\,\pi^*V)^\q$
is a Laurent
polynomial. We see that the natural imbeddings
\[F(V)^\q\hookleftarrow V(\R)^\q\hookrightarrow \Gamma_{hol}(\C^*,\,\pi^*V)^\q=
\Gamma(\E,\,V)\]
are both bijections.  The theorem follows.
$\quad\square$

\vskip 5mm
\pagebreak[3]
{\bf References}
\vskip 2pt

{\footnotesize
[A] $\;\;\;\;$ \parbox[t]{115mm}{
 Atiyah  M.F., {\em Vector bundles over an elliptic curve} Proc. London. 
Math. Soc. (3), 7, 1957, 414-452.}

[B] $\;\;\;\;$ \parbox[t]{115mm}{
Biswas I. {\em Principal bundles admitting a holomorphic connection.}
preprint alg-geom/9601019.}

[BV] $\;\;$ \parbox[t]{115mm}{
Babbit D.G., V.S.Varadarajan. {\em Formal reduction theory of
meromorphic differential equations: a group theoretic view.} Pac.
Jour. Math.
{\bf 109} (1983), no.1, 1-80.}

[DM] $\;\;$ \parbox[t]{115mm}{
Deligne P., Milne J.
{\it Tannakian Categories}, LN in Math, {\bf900} (1982) 101--228.}

[EFK] $\;$ \parbox[t]{115mm}{
Etingof P., Frenkel I., Kirillov A, Jr. {\it Spherical functions
on affine Lie groups.} Duke Math. J., {\bf 80} (1995), 59-90.}

[E] $\;\;\;\;$ \parbox[t]{115mm}{Etingof P.
{\it Galois groups and connection matrices 
for difference equations.} Electronic Research Announcements
{\bf 1} (1995).}

[R] $\;\;\;\;$  \parbox[t]{115mm}{Ramanathan A.
{\em Stable principal bundles on a compact Riemann 
surface.} Math. Ann. {\bf 213} (1975) 129-152.}

[RR] $\;\;$ \parbox[t]{115mm}{
Ramanan S., Ramanathan A. {\em Some remarks on the instability flag.}
Toh\^oku Math. J. {\bf 36} (1984) 269-291.}

[S] $\;\;\;\;\,$ \parbox[t]{115mm}{
Springer T.A. {\em Linear algebraic groups.} Progress in Mathematics vol.9,
Boston: Birkh\"auser, 1981.} 
}

\vskip 5mm
\nab
{\sc University of Chicago, Department of Mathematics, Chicago IL 60637.}
\vskip 4pt\nab
E-mail:\,
{\footnotesize $\;\;
 {barashek@math.uchicago.edu}\quad,\quad {ginzburg@math.uchicago.edu}$}

\end{document}